\newcommand{\Label}[1]{\label{#1}\hspace{.3cm}\fbox{\rm #1}\hspace{.3cm}}
\renewcommand{\Label}{\label} %% COMMENT THIS LINE OUT TO INCLUDE LABELS
\newcommand{\marker}[1]{\mbox{\Huge$\bullet$}#1\mbox{\Huge$\bullet$}}
\renewcommand{\marker}[1]{\mbox{}}
\newcommand{\Comments}[1]{\\ \fbox{\fbox{\rm #1}}\\}
\renewcommand{\Comments}[1]{\mbox{}}
\documentclass[11pt]{article}
 \usepackage{amsmath,amsfonts,amssymb}
 \usepackage{amsfonts}
 \usepackage{amssymb}
\newcommand{\ls}[1]
   {\dimen0=\fontdimen6\the\font \lineskip=#1\dimen0
 \advance\lineskip.5\fontdimen5\the\font \advance\lineskip-\dimen0
 \lineskiplimit=.9\lineskip \baselineskip=\lineskip
 \advance\baselineskip\dimen0 \normallineskip\lineskip
 \normallineskiplimit\lineskiplimit \normalbaselineskip\baselineskip
 \ignorespaces }
 \numberwithin{equation}{section}
 \newtheorem{theorem}{Theorem}[section]
 
 \newtheorem{corollary}[theorem]{Corollary}
 \newtheorem{lemma}[theorem]{Lemma}
 
 \newtheorem{proposition}[theorem]{Proposition}
 \newtheorem{definition}[theorem]{Definition}
 \newtheorem{definitions}[theorem]{Definitions}
 
 \newtheorem{remarks}[theorem]{Remarks}

 \newcommand{\beq}{\begin{equation}}
 \newcommand{\eeq}{\end{equation}}
 \newcommand{\beqa}{\begin{eqnarray}}
 \newcommand{\eeqa}{\end{eqnarray}}
 \newcommand{\beqann}{\begin{eqnarray*}}
 \newcommand{\eeqann}{\end{eqnarray*}}

 \newcommand{\proof}{\noindent {\bf Proof:\ }}
 \newcommand{\pf}{\noindent \mbox{{\bf Proof}: }}
 \def\squarebox#1{\hbox to #1{\hfill\vbox to #1{\vfill}}}
 \newcommand{\qed}{\hspace*{\fill}
    \vbox{\hrule\hbox{\vrule\squarebox{.667em}\vrule}\hrule}\smallskip}
 \newcommand{\req}[1]{(\ref{#1})}
 \newcommand{\lip}{\langle}
 \newcommand{\rip}{\rangle}

 \newcommand{\won}{{\boldsymbol 1}}

 \newcommand{\grad}[1]{\nabla #1}

 \newcommand{\half}{{\frac{1}{2}}}%

 \newcommand{\dfr}{\displaystyle\frac}
 \newcommand{\hilbert}{\bigcirc\kern -0.8em{\rm\scriptstyle {H}\;}}
 \newcommand{\al}{\alpha}
 
 \newcommand{\eps}{\varepsilon}
 \newcommand{\om}{\omega}
 \newcommand{\Om}{\Omega}
 \newcommand{\vph}{\varphi}

\newcommand{\boldR}{\mathbf{R}}
 \newcommand{\DD}{\mathbb{D}}

 \newcommand{\FF}{\calF}
 \newcommand{\FFF}{F}
 
 \newcommand{\calF}{{\mathcal F}}%
 \newcommand{\calL}{{\mathcal L}}%
 \newcommand{\PPi}{\Pi\hspace{-.31cm}\Pi}
 \newcommand{\PPI}{\PPi}
 \newcommand{\calP}{{\mathcal P}}%
 \newcommand{\calS}{{\mathcal S}}%
 \newcommand{\Nat}{\mathbb{N}}%{\mbox{${\rm I\!N}$}}
 \newcommand{\Reals}{\mathbb{R}}%
 \newcommand{\ip}[2]{\left(#1,\rule[0.1cm]{0cm}{0.24cm}#2\right)}
 \newcommand{\pair}[4]{\raisebox{-1.2ex}{\mbox{\tiny$#1$}}\!
       \langle#2,#3\rangle\!\raisebox{-1.2ex}{\mbox{\tiny $#4$}}}
 \newcommand{\bpair}[4]{\raisebox{-2.4ex}{\mbox{\tiny$#1$}}\!
       \left\langle#2,#3\right\rangle\!\raisebox{-2.4ex}{\mbox{\tiny $#4$}}}
 
 \newcommand{\clfs}[2]{\pair{W^*}{#2}{#1}{\Wss}}

 \newcommand{\dom}{{\rm dom}}
 \newcommand{\ddom}{{\rm \bf{dom}}}
 
 \newcommand{\domdel}[1]{\dom_{#1}\delta}
 \newcommand{\ddomdel}[1]{\ddom_{#1\,}\ddelta}
 \newcommand{\Domdel}[1]{{\ddomdel{#1}}}
 \newcommand{\ddelta}{\delta}

 \newcommand{\Dp}{\mathbb{D}_{p,1}}

 \newcommand{\Lwa}[1]{L^2_{\rm wa}\!\left(\mu;\rule[0.1cm]{0cm}{0.25cm}L(H,#1)\right)}
 \newcommand{\Wss}{W^{**}}
 
 \newcommand{\DDD}{\mathbf D}
 \newcommand{\KKK}{\mathbf K}
 \newcommand{\QQQ}{\mathbf Q}
 \newcommand{\RRR}{\mathbf R}
 \newcommand{\TTT}{\mathbf T}

 \newcommand{\tr}{{\rm tr\,}}
 \newcommand{\dbr}[2]{\langle\!\langle#1,#2\rangle\!\rangle}
\begin{document}
 \thispagestyle{empty}
 \noindent
% {\em Erratum}\\
  {\em The original paper follows this $1$-page erratum}
  \bigskip

 \noindent
 {\bf \large Erratum to ``The Clark-Ocone formula for vector valued random variables
        in abstract Wiener space", Jour. Func. Anal. 229, 143--154  (2005)}

 \noindent {\small  \bf E. Mayer-Wolf\footnote{Department of Mathematics, Technion, Israel;\ \ emw@tx.technion.ac.il},
                 M. Zakai\footnote{Department of Electrical Engineering, Technion, Israel;\ \ zakai@ee.technion.ac.il}}
 \vspace{.5cm}

  In this paper we considered the extension of the Clark-Ocone formula for a
  random variable defined on an abstract Wiener space $(W,H,\mu)$ and taking
  values in a Banach space (denoted there either $B$ or $Y$). The main result
  appears in Theorem~3.4. Unfortunately, as first pointed out to us by J. Maas
  and J. Van Neerven, the dual predictable projection $\PPI$ introduced in
  Definition~3.1(iii) via the characterization~(3.1), does \underline{not}
  define a random operator in $L^2(\mu;L(H,Y))$ as claimed, but rather an
  element of the larger space $L(H,L^2(\mu,Y))$. Consequently the right hand
  side of~(3.6) in the main result is ill defined.

  We have been unable to overcome this difficulty in a meaningful way. We
  should point out that a Clark-Ocone formula was recently obtained
  in~\cite{MvN} for random variables on a classical cylindrical Wiener space
  taking values in a UMD Banach space, in which $\delta$ can be explicitly
  defined \`{a} la It\^{o} on adapted processes. Our work, however, was
  different in spirit and made use of the extended version of $\delta$
  introduced in \cite{div}. While it is possible to provide an even weaker
  interpretation of~(3.6) in which $\delta$ is extended to suitable elements
  of $L(H,L^2(\mu;Y))$, the result would have amounted to little more than the
  collection of classical Clark-Ocone formulae for the scalar random variables
  $\{\langle v,y^*\rangle,\,y^*\!\in\!Y^*\}$.

  The main result,~Theorem 3.4, is thus considerably weakened; it remains true
  a) assuming that $Y^{**}$ has the Radon Nikodym property (RNP) with respect
  to $\mu$, and b) for $Y$\!-\!valued random variables $v$ for which one can
  verify that $\PPI\grad v\!\in\!L(H,L^2(\mu;Y))$.\ (The need for the
  additional RNP condition a) derives from an error, also brought to our
  attention by J. Maas and J. Van Neerven, in the proof of Proposition~3.14
  of~\cite{div}, (cited here as Lemma~2.3) which has been corrected
  in~\cite{err} under the RNP condition).

  Section~4 is not affected by the difficulties described above.

 \ls{1.3}
   \title{The Clark--Ocone formula for vector valued random variables in
          abstract Wiener space}
   \author{E. Mayer-Wolf
           \footnote{Department of Mathematics, Technion I.I.T., Haifa, Israel}
       \ \ and M. Zakai\footnote{Department of Electrical Engineering,
                                       Technion I.I.T., Haifa, Israel}}
%   \date{\today}
 \date{}
 \maketitle
 \setcounter{page}{1}
%  \centerline{\huge Draft}
%  \bigskip
%
 \begin{abstract}
  The classical representation of random variables as the It\^{o} integral of
  nonanticipative integrands is extended to include Banach space valued random
  variables on an abstract Wiener space equipped with a filtration induced by
  a resolution of the identity on the Cameron--Martin space.
  The It\^{o} integral is replaced in this case by an extension of the
  divergence to random operators, and the operators involved in the
  representation are adapted with respect to this filtration in a suitably
  defined sense.
 \end{abstract}
%
%\vspace{1cm}
 \vfill \noindent{\sc Key Words and Phrases: Clark formula, Clark-Ocone
     formula, Banach space valued random variables, weakly adapted operators, }
     % Abstract Wiener Space, Divergence, Flows.\\
 %; Secondary 60H05.\\

\newpage
 \section{Introduction}

  The representation of square integrable functionals of the Wiener process
  as a sum of multiple Wiener-It\^{o} integrals was derived by K. It\^{o} in
  his 1951 paper \cite{ito}. It follows easily from this series that
  every such functional is representable as a Ito\^{o} integral.
  This representation, however, was not stated explicitly in \cite{ito}, and
  its first appearance seems to have occurred in the 1967 paper of H. Kunita
  and S. Watanabe \cite{kuwa}.

  The problem of finding an explicit expression for the integrand in the
  It\^{o} integral was formulated and solved under certain differentiability
  restrictions by J. M. C. Clark in 1970 \cite{clark}. In 1984, D. Ocone
  \cite{oco} applied the Malliavin calculus to relax these restrictions
  significantly, and then in further generality with I. Karatzas and J. Li
  \cite{kol}. In loose terms, this representation is valid for $L^2$ (more
  generally, $L^1$) random variables\ $\vph$ on Brownian paths\
  $\om\!=\!(\om_t)_{0\le t\le 1}$, smooth enough that there exists a
  (``derivative") process  $D_t\vph$ such that
  \[ \left.\frac{d\vph\left(\om\!+\!\eps\int_0^\cdot h_s\,ds\right)}{d\eps}
                \right|_{_{\eps=0}} =\int_0^1 D_t\vph\,h_t\,dt \]
  in an appropriate sense. The Clark--Ocone formula then states that
  \[ \vph=E\vph+\int_0^1 E\left(\left.D_t\vph\right|\FF_t\right)\,d\om_t,\]
  where\ $(\FF_t)$ is the canonical filtration.

  The purpose of this paper is to obtain the Clark representation for random
  variables taking values in Banach spaces. This will be done in the
  context of an abstract Wiener space\ $(W,H,\mu)$\ whose natural filtration is
  induced by a resolution of the identity, thus allowing for the notion of
  adaptedness. Extensions of the Clark--Ocone formula in an abstract Wiener
  space have already been studied (\cite{wuli},\cite{uz1},\cite{oss}) from a
  different
  point of view, namely, for scalar random variables.

%  Here add an outline...
%
  Section~2 is devoted to some basic notions of stochastic analysis in Wiener
  space, including the gradient and divergence operators, the latter applied
  to random variables which are not necessarily $H$-valued, as introduced
  in~\cite{mwz}.
  In Section~3 we first summarize the necessary preliminaries concerning
  resolutions of the identity, their induced filtrations and vector valued
  random variables adapted with respect to them, based mostly on \cite{uz1},
  \cite{uz2} and~\cite{zak}.
  Next we consider the divergence of (weakly adapted) random variables taking
  values in a Banach space $B$ (which reduces to the It\^{o} integral when $B$
  is the Cameron  Martin space) and then apply these results and those of
  Section~2 to derive %first the It\^{o} representation of second order
  %Banach space valued variables, and then
  the Clark--Ocone formula for those
  such variables which are regular . This will be illustrated in Section~4 where
  measure preserving transformations on Wiener space are considered as $W$-valued
  random variables. Section~5 contains some concluding remarks.
  \section{Stochastic analysis preliminaries}
  An abstract Wiener space $(W,H,\mu)$ consists of a separable Banach space
  $W$, a separable Hilbert space $H$ densely embedded in $W$ and a zero mean
  Gaussian measure $\mu$ on $W$'s Borel sets under which each $l\!\in\!W^*$
  is a ${\rm N}(0,|l|_{\mbox{\tiny H}}^2)$ random variable, denoted $\delta l$.
  Here $W^*$ was implicitly taken to be a dense subspace of\ $H$, as it will be
  throughout. By density, this extends to a zero mean linear Gaussian random
  field\ $\{\delta h,\,h\!\in\!H\}$ whose covariance is induced by\ $H$'s
  inner product.

  Let $(\eta_n)$ be an independent sequence of $N(0,1)$ random variables on
  some probability space $(\Om,\FF,P)$, and $(e_n)$ an orthonormal base (ONB)
  of\ $H$. It\^{o}-Nisio's theorem \cite{itni} states that\
  $\sum_{n=1}^\infty \eta_ne_n$ converges to a $W$-valued random variable $\xi$
  whose distribution is $\mu$, and that if in particular\  $\Om\!=\!W$\ and\
  $\eta_n\!=\!\delta e_n$\ for each $n$,\ then $\xi(w)\!=\!w\ \ \mu$ a.s.

  For any Banach space $Y$ and $1\!\le\!p\!\le\!\infty$\ we denote by
  \ $L^p(\mu;Y)$ the class of strongly measurable $Y$-valued random variables
  $v$ on $W$ such that $\|v\|_{\mbox{\tiny Y}}\!\in\!L^p(\mu)$, and
   \vspace{-.8cm}

 \begin{equation} \Label{Sn}
  \calS(Y)=\mbox{\Large {$\{$}}F\!%=\sum_{j=1}^m \Phi_jb_j:
    :=\!\sum_{j=1}^m\underbrace{\vph_j(\delta h_1,\ldots,\delta h_n)}_{\Phi_j}b_j\
      \mbox{\Large {$|$}}\ m,n\!\in\!\Nat,\ \vph_j\!\in\!C^\infty_b(\Reals^n),
      \ h_i\!\in\!H,\ b_j\!\in\!Y\ \mbox{\Large {$\}$}},
 \end{equation}
 \vspace{-.8cm}

  and the gradient of these {\em simple} $Y$-valued random variables is defined
  to be\vspace{-.5cm}

 \begin{equation} \Label{simplegradient}
  \grad F =\sum_{j=1}^m \grad\Phi_j\otimes b_j=\sum_{j=1}^m\sum_{i=1}^n
       \partial_i \vph_j(\delta h_1,\ldots,\delta h_n)\,h_i\otimes b_j
       \ \ \in L^\infty(\mu;L(H,Y)).
 \end{equation}
  Here and throughout $L(X,Y)$ denotes the space of bounded linear
  operators from a Banach space $X$ to a Banach space $Y$, equipped with their
  operator norm\ (and $L(X)\!=\!L(X,X)$). It should be noted that when $Y$
  is a separable Hilbert space, the Hilbert--Schmidt norm of $\grad F$ is
  traditionally used; the operator norm in this case was first considered by
  G. Peters in~\cite{pet}.%, even when $Y$ is a separable Hilbert space.

%  Now let $Z\!\subset\!W$\ be a linear subspace (\,$L(W,Y)\!\subset\!L(Z,Y)$
%  naturally),
   For each \ $1\!\le\!p\!<\!\infty$ define on\ $\calS(Y)$ the norms
   \begin{equation}
   \|\FFF\|_{p,1}=  \left(\|\FFF\|^p_{L^p(\mu;Y)}
              +\|\grad{\FFF}\|^p_{L^p(\mu;L(H,Y))}\right)^{\frac{1}{p}}.
   \end{equation}
  The Sobolev spaces $\Dp(Y)\!\subset\!L^p(\mu;Y)$\ are defined to be\
  $\calS(Y)$'s completions according to these norms. By closability, $\grad$
  can be extended to a bounded operator (with a slight abuse of notation)
  $\grad\!:\!\Dp(Y)\to L^p(\mu;L(H,Y))$.
%  The subspace $Z$ of ``differentiability directions" should be thought of
%  as $W$ when imposing Frech\'{e}t differentiability (which is usually too
%  restrictive), as $H$ in the classical Malliavin calculus, or as $W^*$ when
%  even weaker smoothness is imposed, as needed when considering the divergence
%  of $W$-valued random variables, for example (c.f. below and \cite{mwz}).

  The divergence operator on random operators in\ $L(H,Y)$ is defined by
  duality. Recall that the trace \ ${\rm tr}\,\TTT$\ of an operator\
  $\TTT\!\in\!L(H)$, which is defined to be\
  $\sum_{n=1}^\infty \clfs{\TTT e_i}{e_i}$ if this sum converges and is the
  same for every ONB $(e_n)$ of $H$, induces  the pairing
  $\dbr{\KKK}{\DDD}\!:=\!\tr\left(\KKK^T\DDD\right)$, for appropriate\
  $\KKK\!\in\!L(H,Y)$ and $\DDD\!\in\!L(H,Y^*)$.
% \ (in particular if at least one of the operators $\KKK,\DDD$ has finite
% rank or, more generally, is nuclear).
%
  We shall say that\ $\KKK\!\in\!L^1(\mu;L(H,Y))$ has finite rank
  if for some\ $m\!\in\!\Nat$,\ \ $\KKK\!=\!\sum_{k=1}^mh_j\otimes
  y_j$\ with\ $u_j\!\in\!L^1(\mu;H)$\ and\ $y_j\!\in\!Y$,\ that is,\
  $\KKK h\!=\!\sum_{j=1}^m (u_j,h)y_j$.% If either $\KKK$\ or\ $\DDD$\ has
%  finite rank then $\dbr{\KKK}{\DDD}$\ is well defined.
  \begin{definition} \Label{Domp}
   \hspace*{.2cm}For $1\!\le\!p\!<\!\infty$\ \ \ let\ \
   $\Domdel{p,Y}$%$\!=\!\Domdel{p}\!$
   \ \ be the set of all\ \ $\KKK\!\in\!L^p(\mu;L(H,Y))$\ \ for which there
   exists a\ \ $\ddelta \KKK\!\in\!L^p(\mu;Y^{**})$,\ \ the {\bf divergence} of\
   $\KKK$,\ \ such that for all\ $\FFF\!\in\!\calS(Y^*)$.
   \begin{equation} \Label{opdual}
    E\,\dbr{\KKK}{\grad{\!\FFF}}=E\pair{Y^*}{F}{\delta \KKK}{Y^{**}}
   \end{equation}
  \end{definition}
%    on the $H$-valued ones among them, $\delta$ coincides with the usual divergence.
  (Note that the pairing in \req{Domp} is well defined since\ $\grad{\!\FFF}$
  has finite rank). A necessary and sufficient condition for
  $\KKK\!\in\!\Domdel{p,Y}$ (cf. \cite[Equation (3.12)]{mwz}) is
  that for some $\gamma\!>\!0$
  \[|E\,\dbr{\KKK}{\grad{\!\FFF}}|\!\le\!\gamma\|\FFF\|_{_{L^q(\mu;Y^*)}} \]
  ($\frac{1}{p}\!+\!\frac{1}{q}\!=\!1$)\ for all $\FFF\!\in\!\calS(Y^*)$.

  \noindent
  Lemma~\ref{weakversion} below provides a ``weak"
  characterization of $\delta\KKK$. If $\delta$ had been required to be
  $Y$--valued (and not only $Y^{**}$--valued), the ``if" implication in the
  Lemma would no longer be valid.
  \vspace{-.1cm}

  \noindent
   We denote\ $\Domdel{p,\Reals}\!=\!\domdel{p}$; this space contains $H$-valued
  random variables, and in this case $\delta$ is the usual divergence.
  \begin{remarks}\Label{delten}\mbox{}\vspace{-.5cm}

   \begin{itemize}
    \item[i)]\cite[Remark 3.13]{mwz}\ If $\KKK$'s range is $\mu$-a.s. contained
     in a (deterministic) finite dimensional subspace of $Y$,\ \req{opdual}
     extends to all\ $F\!\in\!\Dp(Y^*)$.\vspace{-.3cm}

    \item[ii)] If $\alpha\!\in\!\domdel{p}$\ and $y\!\in\!Y$, it follows
     directly from the definitions that\ \
     $\alpha\otimes y\!\in\!\Domdel{p,Y}$\ \
     and that $\delta(\alpha\otimes y)\!=\!(\delta\alpha)y$.
   \end{itemize}
  \end{remarks}

  \begin{lemma}\cite[Proposition 3.14]{mwz}\Label{weakversion}
    An element\ $\KKK\!\in\!L^p\left(\mu;L(W^*,Y)\right)$\ belongs to
    $\Domdel{p,Y\,}$\ \ if and only if \ \  $\KKK^Tl\in\domdel{p}$\ \
    for every $l\!\in\!Y^*$\ and for some\ $C\!>\!0$
    \begin{equation} \Label{Cbound}
     \|\ddelta \left(\KKK^Tl\right)\|_{_{L^p(\mu)}}\!\le\!C\|l\|_{_{Y^*}}
                                        \hspace{1.5cm}\forall\ l\!\in\!Y^*.
    \end{equation}
    In this case
%\begin{subequations}
   \begin{equation} \Label{weak}
    \delta(\KKK^Tl)=\pair{Y^*}{l}{\ddelta \KKK}{Y^{**}}\hspace{1.5cm}{\rm a.s.}
   \end{equation}
  and more generally, for any\ $F\!\in\!\calS(Y^*),\ \ \ \
  \KKK^TF\!\in\domdel{p}$\ \ \ and
 \begin{equation} \label{weakb}
  \delta(\KKK^TF)
      =\pair{Y^*}{F}{\ddelta \KKK}{Y^{**}}-\dbr{\KKK}{\grad^{^{\!W^*}}\!\!F}.
 \end{equation}
%\end{subequations}
  \end{lemma}
 \noindent{\bf Examples}\vspace{-.5cm}

% \small
 \begin{itemize}
  \item[i)] If $v(w)\!\equiv w_0\!\in\!W$\ belongs to $\domdel{1}$,\ then
   necessarily $w_0\!\in\!H$\ \cite[Remark 3.2b)]{mwz}.
  \item[ii)] $v(w)\!=\!w$ does not belong to $\domdel{1}$. This follows by
   applying \cite[Proposition~3.6)]{mwz} to the It\^{o}-Nisio representation
   $v\!=\!\sum_n \delta e_n\,e_n$\ for any ONB $(e_n)$.
  \item[iii)]
   $v(w)\!=\!\sum_{n=1}^\infty (\delta e_{2n}\,e_{2n-1}-\delta e_{2n-1}\,e_{2n})$
   converges and, like in (ii), $v\stackrel{\mathcal D}{\sim}\mu$\ (by the
   It\^{o}-Nisio theorem). However, here $v\!\in\!\domdel{1}$ and $\delta v=0$.
   This follows from \cite[Lemmas~3.3, 3.4]{mwz}.
  \item[iv)] $\won_H$ %(viewed as $W^*\rightarrow \Wss$)
   belongs to $\Domdel{p,W}$ for all $p\!\ge\!1$ (but not to $\Domdel{p,H}$\ !)
   \ and $\ddelta\won_H(w)\!=\!w$\ \ $\mu$-a.s. \cite[Corollary~3.16)]{mwz}.
 \end{itemize}
 \section{Adaptedness and the divergence representation of vector-valued
   random variables}
  Let $\pi\!=\!\{\pi_\theta,\,\theta\!\in\![0,1]\}$\ be a strictly increasing
  continuous resolution of the identity on $H$ (the $\pi_\theta$'s are
  orthogonal projections in $H$ with $\pi_0\!=\!0$\ and\ $\pi_1\!=\!I_H$). Each
  such resolution of the identity induces the filtration\
  $\FF=\{\FF_\theta,\theta\!\in\![0,1]\}$ on\ $(W,H,\mu)$\ defined by
  \[ \FF_\theta=\sigma\left(\rule[0.1cm]{0cm} {0.27cm}\delta(\pi_\theta h),
                     \ h\!\in\!H\right)\hspace{1.5cm}\theta\!\in\![0,1] \]
  which generates a time structure with respect to which notions of
  adaptedness can be defined.\\

  \noindent
  {\bf a.\ \ }\underline{\bf Adaptedness}
  \begin{definitions} \Label{defadapted}
   Let $Y$ be an arbitrary Banach space.\vspace{-.5cm}

   \begin{itemize}
    \item[i)] An $H$-valued random variable $u$ is {\em adapted} (to $\FF$)\ \
     if\ $\left(u,\rule[0.1cm]{0cm}{0.24cm}\pi_\theta h\right)$ is
     $\FF_\theta$--measurable for each $h\!\in\!H$ and $\theta\!\in\![0,1]$.
     \ \ \ Set\ \ $L^2_{\rm a}(\mu;H)\!=\!
     \left\{u\!\in\!L^2(\mu;H),\rule[0.1cm]{0cm}{0.27cm}\
                                        u\ {\rm is\ adapted}\right\}$.
    \vspace{-.35cm}

    \item[ii)] An $L(H,Y)$--valued random operator $G$ is {\em weakly adapted}
     (to $\FF$)\ \ if\ \ $G^Ty^*$ is adapted to $\FF$ for  each $y^*\!\in\!Y^*$.
      \ \ \ Set\ \ $\Lwa{Y}\!=\!
        \left\{G\!\in\!L^2\left(\mu;\rule[0.1cm]{0cm}{0.25cm}L(H,Y)\right),\
             \rule[0.1cm]{0cm}{0.3cm}G\ {\rm is\ weakly\ adapted}\right\}$.
    \vspace{-.35cm}

    \item[iii)]
     $\Pi=$ the orthogonal projection of $L^2(\mu;H)$\ onto\
     $L^2_{\rm a}(\mu;H)$\ \ \ \ and\\
     \marker{change\ 1\ -\ 1/7}\\
     $\PPi\!:\!L^2\left(\mu;\rule[0.1cm]{0cm}{0.25cm}L(H,Y)\right)
           \longrightarrow\Lwa{Y}$\ is defined by
     \begin{equation} \Label{Pidef}
      \pair{Y^*}{y^*}{(\PPi\KKK) h}{Y}
         =\left(\rule[0.1cm]{0cm}{0.26cm}\Pi(\KKK^Ty^*),h\right)_H,\hspace{1cm}
           \KKK\!\in\!L^2\left(\mu;\rule[0.1cm]{0cm}{0.25cm}L(H,Y)\right),\
           h\!\in\!H,\ y^*\!\in\!Y^*.
     \end{equation}
   \end{itemize}
  \end{definitions}
%  \begin{remark}
  It follows directly from \req{Pidef} that
  \begin{equation} \Label{weakproj}
   \Pi\!\left(\KKK^Ty^*\right)\!=\!\left(\PPi\KKK\right)^T\!y^*\hspace{1.5cm}
     \forall\KKK\!\in\!L^2\left(\mu;\rule[0.1cm]{0cm}{0.25cm}L(H,Y)\right),
     y^*\!\in\!Y^*
  \end{equation}
   from which it follows that $\PPi\KKK$ is indeed weakly adapted for every
   $\KKK\!\in\!L^2\left(\mu;\rule[0.1cm]{0cm}{0.25cm}L(H,Y)\right)$. %Moreover,
   $\PPi$ is a {\em projection} onto $\Lwa{Y}$,\ as can be easily verified,
   which moreover inherits from $\Pi$ the weak orthogonality property
   \begin{equation} \Label{weakortho}
    E\dbr{\KKK}{\QQQ}=E\dbr{\PPi \KKK}{\QQQ}
   \end{equation}
   for every $\KKK\!\in L^2(\mu;L(H,Y)$\ and finite rank\ $\QQQ\!\in\!\Lwa{Y^*}$.
   Indeed, if\ $Q\!=\!q\otimes y^*$,\ with\  $q\!\in\!L^2_a(\mu;H)$\ and\
   $y^*\!\in\!Y^*$,\ then
   \[ E\tr{\KKK^T(q\otimes y^*)}=E\tr{q\otimes\KKK^Ty^*}
                                         =E(q,K^Ty^*)=E(q,\Pi K^Ty^*), \]
   since\ $q$\ is adapted, and the same expression is obtained
   when $\KKK$ is replaced by\ $\Pi\KKK$.
   \ \ \ \marker{1}
%  \end{remark}

  The following lemma suitably generalizes the It\^{o} integral of adapted
  processes, and its isometry property
  \begin{equation} \Label{Itoisometry}
   \marker{change\ 2 - 1/7}\ \ \
   E\delta(u)\delta(v)\!=\!E(u,v)\hspace{1cm}
      \forall u,v\!\in\!L^2_a\left(\mu;\rule[0.1cm]{0cm}{0.25cm}L(H,Y)\right).
   \ \ \ \marker{2}
  \end{equation}
  A random operator $G(\omega):X\to Y$ has finite rank if
  $G\!=\!\sum_{j=1}^m x^*_j\otimes y_j$ for appropriate\
  $m\!\in\!\Nat$,\ $X^*$-valued random variables $x^*_j(\omega)$\ and\
  nonrandom $y_j\!\in Y$,\ $1\!\le\!j\!\le\!m$.
  \begin{lemma} \Label{itoisometry}
   \mbox{}\vspace{-.4cm}

   \begin{itemize}
     \item[i)] For any Banach space $Y$,\ \ \ $\Lwa{Y}\!\subset\Domdel{_{2,Y}}$.
     If, moreover, $\DDD\!\in\!\Lwa{Y}$ has finite rank, then
     $\delta\DDD\!\in\!Y$. \vspace{-.3cm}

    \item[ii)] Given a Banach space $B$, if $\KKK\!\in\Lwa{B}$\ \ \ and\ \ \
        $\DDD\!\in\!\Lwa{B^*}$\ has finite rank, then
    \begin{equation} \Label{isometry}
     E\pair{B*}{\delta\DDD}{\delta\KKK}{B^{**}}=E\dbr{\KKK}{\DDD}.
    \end{equation}
   \end{itemize}
  \end{lemma}
%  (\small The tensor notation for $\DDD$ means that
%   $\DDD h=\sum_{j=1}^n\!\left(\psi_j,h\right)\,y_j^*$\ \normalsize).
%
  \pf
   For any $G\!\in\!\Lwa{Y}$ and\ $y^*\!\in\!Y^*$, it holds by definition that
   $G^Ty^*\!\in\!L^2_{\rm a}(\mu;H)$. It is well known that adapted $H$--valued
   random variables of second order are It\^{o} integrable, and thus in
   $\domdel{2}$. Lemma~\ref{weakversion} then implies that
   $G\!\in\!\Domdel{_{2,Y}}$.\vspace{-.2cm}

   \noindent
   If $\DDD\!=\!\sum_{j=1}^m\!\vph_j\otimes y_j$, with
   $\vph_j\!\in\!L^2_{\rm a}(\mu,H)$\ \ and\
   $y_j\!\in\!Y$ ,\ $1\!\le\!j\!\le\!m$,\ \ then by Remark~\ref{delten}\,ii)
   $\DDD\!\in\!\Domdel{_{2,Y}}$\ and\
   $\delta\DDD\!=\!\sum_{j=1}^m(\delta\vph_j)y_j$.

   \noindent
   As for ii), let\ $\DDD\!=\!\sum_{j=1}^m u_j\otimes b^*_j$,\ with
   $u_j\!\in\!L_{\rm a}(\mu;H)$ and\ $b_j^*\!\in\!B^*$,\ $1\!\le\!j\!\le\!m$,
   and let $(e_i)_{i\in\Nat}$ be an arbitrary ONB in $H$. Then
   \begin{eqnarray*}
    E\pair{B*}{\delta\DDD}{\delta\KKK}{B^{**}}
     &=&E\sum_{j=1}^m\delta u_j\pair{B^*}{b_j^*}{\delta\KKK}{B^{**}}\\
     \marker{change\ 3 - 1/7}
     &\stackrel{\req{weak}}{=}&\ \sum_{j=1}^mE\delta(u_j)\,\delta(\KKK^Tb^*_j)\\
     &\stackrel{\req{Itoisometry}}{=}&\marker{3}\sum_{j=1}^mE\ip{u_j}{\KKK^Tb_j^*}
        =E\sum_{j=1}^m\sum_{i=1}^\infty\ip{u_j}{e_i}\ip{e_i}{\KKK^Tb^*_j}\\
     &=&E\sum_{j=1}^m\sum_{i=1}^\infty\ip{u_j}{e_i}\pair{B}{\KKK e_i}{b^*_j}{B^*}
        =E\sum_{i=1}^\infty\bpair{B}{\KKK e_i}{\sum_{j=1}^m\ip{u_j}{e_i}{b^*_j}}{B^*}\\
     &=&E\sum_{i=1}^\infty \pair{B}{\KKK e_i}{\DDD e_i}{B^*}=E\dbr{\KKK}{\DDD}.
   \end{eqnarray*}
   \nopagebreak
   \vspace{-1.72cm}

   \qed
   \begin{corollary} \Label{divfree}
    If\ $\KKK\!\in\!\Lwa{Y}$\ and\ $\delta\KKK\!=\!0$\ \ \ then\ \ \
    $\KKK\!=\!{\mathbf 0}$.
   \end{corollary}
   \pf
    Under the assumptions on $\KKK$ it follows from~\req{isometry} that
    $E\dbr{\KKK}{\DDD}\!=\!0$\ for every finite range weakly adapted
    random operator $\DDD\!:\!H\!\to\!B^*$, in particular\
    $\DDD\!=\!\vph\!\otimes\!b^*$\ with\ $\vph\!\in\!L^2_{\rm a}(\mu;H)$\
    and $b^*\!\in\!B^*$. Thus
    \[ 0=E\dbr{\KKK}{\DDD}=E\ip{\vph}{\KKK^Tb^*},\]
    and since $\vph,\ b^*$ were arbitrary, the conclusion follows.
   \qed
   \vspace{.5cm}

   \noindent
   \noindent
   {\bf b.\ \ }\underline{\bf The Clark--Ocone formula}\\

   \noindent
%   When the random variable $v$ in Theorem~\ref{Itorepresentation}
%   is sufficiently regular, the random operator $\KKK$ in its
%  representation\ can be written somewhat more explicitly.
   This subsection is devoted to the main result of this note.
   \begin{theorem}  \Label{clarkocone}
    Given a Banach space $B$\ and\ $v\!\in\!{\mathbb D}_{2,1}^H(B)$,
    \begin{equation} \Label{clarkformula}
     v=Ev+\delta\left(\PPi\grad v\right)
    \end{equation}
    and $\KKK\!=\!\PPi\grad v$ is the unique element in $\Lwa{B}$ such that
    $v\!=\!Ev\!+\!\delta\KKK$.
   \end{theorem}
   (By Lemma~\ref{itoisometry}i),\ $\PPi\grad v$\ indeed belongs to\
    $\delta$'s\ domain.)

   \pf
    We shall again assume that $Ev\!=\!0$.
    Let $F\!=\!\sum_{i=1}^n \Phi_ib_i^*\in \calS(B^*)$ be a simple random
    variable (c.f.~\req{Sn}) for which $E\Phi_i\!=\!0$ for each $i$. By the
    standard It\^{o} representation, $\Phi_i\!=\!\delta q_i$,\ for
    appropriate\ $q_i\!\in\!L^2_{\rm a}(\mu;H)$,\ \ $i\!=\!1,\ldots,n$, so that
    \[ F=\sum_{i=1}^m\delta(q_i)\,b^*_i=\delta(\QQQ)\hspace{1.5cm}{\rm with}
            \ \ \ \  \QQQ=\sum_{i=1}^m q_i\otimes b^*_i\ \in\Lwa{B^*}.\]
    We shall show that\
    \begin{equation} \Label{toprove}
     \pair{B}{v}{F}{B^*}\!=\!E\pair{B}{\delta\left(\PPi\grad v\right)}{F}{B^*}
    \end{equation}
    from which \req{clarkformula} will follow since these test variables $F$
    are dense in $L^2(\mu;B^*)$. We have
     \begin{eqnarray*}
      E\pair{B}{v}{F}{B^*}
        &=&E\pair{B}{v}{\delta\QQQ}{B^*}\\
        &=&E\dbr{\grad v}{\QQQ}\\
%     \end{eqnarray*}
%     Denoting $\KKK\!=\!\grad v$, and assuming for simplicity that
%     $\QQQ\!=\!q\!\otimes b^*$\ ,the term on the right is
%     \begin{eqnarray*}
%      E\trace
        &=&E\dbr{\PPi\grad v}{\QQQ}\\
        &=&E\pair{B}{\delta\left(\PPi\grad v\right)}{\delta\QQQ}{B^*}
         = E\pair{B}{\delta\left(\PPi\grad v\right)}{F}{B^*},
     \end{eqnarray*}
     where Remark~\ref{delten}\,i) was used in the second equality,
     \req{weakortho} in the third and Lemma~\ref{itoisometry}\,ii)
     in the fourth.
%     Since $F$ was an arbitrary  element in\ $\calS(B^*)$ with zero
%     expectation, \req{clarkformula} follows.

    \noindent
%    The uniqueness of $\KKK$ has already been established in
%    Theorem~\ref{Itorepresentation}
     As for the uniqueness, if\ $v\!=\!\delta\KKK_i$\ with
     $\KKK_i\!\in\!\Lwa{B},\ i\!=\!1,2$, it follows that
     $\delta\left(\KKK_1\!-\!\KKK_2\right)\!=\!0$\ and thus\ \
     $\KKK_1\!=\!\KKK_2$\ \ by Corollary~\ref{divfree}.
    \qed
 %%%%%%%%%%%%%%%%%%%%%%%%%%%%%%%%%%%%%%%%%%%%%%%%%%%%%%%%%%%%%%%%%%%%%%%%%%%%
%SECTION 4

%\setcounter{section}{3}

 \section{Measure preserving transformations on the Wiener space}

 Let $(W,H,\mu)$ be an abstract Wiener space and let $e_i, i=1,2,\dots$ take
 values in $W^*$ and such that the images of the $e_i$ in $H$ are a
 complete orthonormal base on $H$. By the Ito-Nisio theorem \cite{itni}
 \begin{equation} \label{4-1}
  w_n =\sum_1^n\, \delta (e_i) \, e_i
 \end{equation}
 with $ e_i $ considered as elements in $W$, converges in $ L_1 $ on $W$ to $w$,
 similarly if $ \{\eta_i , i=1,2,\dots\} $ are i.i.d., $ N(0,1) $ then
 $ \sum_1^n\, \eta_i e_i $ converges in $L_1 $ on $W$ to a $W$-valued random
 variable which has the same probability law as $w$. In this case
 $ Tw := \sum_1^\infty\, \eta_i e_i $ will be denoted an ``abstract Wiener
 process'' or ``a measure preserving transformation on the Wiener space'' or
 (for reasons that will become clear later) ``a rotation''.
 Note that $w$ and $ Tw $, while each being Gaussian are, in general, not
 jointly Gaussian. The fact that $ Tw $ as defined above is $W$-valued suggests
 the problem of the Clark representation of this transformation. We have already
 noted that for $ Tw=w , w=\delta (I) $.
 The analysis and characterization of measure preserving transformations is not
 new (\cite{uzbook},\cite{zak}) and most of the results presented here are known;
 it is, however,  more natural to analyze the class of measure preserving
 transformations in the context of this section.

 We prepare the following result for later reference:
 \begin{proposition} \label{prop-4}
  Let $\boldR(w)$ be an a.s. bounded operator on $H$. Assume that $\boldR(w)$ is
  weakly adapted with respect to a filtration induced by a continuous increasing
  $\pi$. Since $\boldR h$ is adapted it is in the domain of $\delta$. Assume
  that the probability law of $\delta(\boldR h)$ is $N(0,|h|_H^2)$, then:
  \begin{enumerate}
   \item If $h_1, h_2\!\in\!H$ and $(h_1,h_2)_H\!=\!0$\ \  then\ \
    $\delta(\boldR h_1)$ and $\delta(\boldR h_2)$  are independent.
   \item $\boldR(w) $ is a.s. an isometry on $H$.
   \item $\sum_i\delta(\boldR e_i)\,e_i$ is measure preserving, and if
    $(e_i)$ and $(h_i)$ are ONB's of $H$ then, a.s.,
    \begin{equation} \label{4-2}
     \sum_i\delta(\boldR h_i)h_i\!=\!\sum_i\delta(\boldR e_i)\,e_i \, .
    \end{equation}
  \end{enumerate}
 \end{proposition}
 \newpage
 \proof

 \hspace*{-.75cm}$ \begin{array}[h]{lll}
 1. & \!\! E\,\exp \{i\al \delta (\boldR h_1)\} \exp \{i\beta\delta (\boldR h_2)\}
    & = E\exp \left\{\delta \Bigl(\al h_1 + \beta h_2)\Bigr)\right\} \\[.2cm]
   && = E\exp \left\{-\dfr{\al^2}{2} \, |h_1|_H^2 -\dfr{\beta^2}{2}\,
       |h_2|_H^2\right\}\\
   [.2cm] && = E\exp\{i\al \delta (\boldR h_1)\} \,
               E\exp \{i\beta\delta (\boldR h_2)\} \, .
\end{array}$

\noindent
2. By part 1, $ y_\theta =\delta (\boldR\pi_\theta h) $ is a Gaussian process of independent increments.
\newline
Hence it is Gaussian martingale and its quadratic variation satisfies
\begin{equation}
\label{4-3}
\lip y,y\rip_\theta = E y_\theta^2 \, .
\end{equation}
and by our assumption $ Ey_\theta^2 = |\pi_\theta h|_H^2 $. But
\begin{equation}
\label{4-4}
\lip y,y\rip_\theta = (\boldR\pi_\theta h , \boldR\pi_\theta h)_H
\end{equation}
and $ \boldR^T \boldR=I $ follows.
\newline
3. Follows from the Ito-Nisio theorem.
\qed

 \begin{theorem} \Label{the-4-1}
  Let $(W,H,\mu)$ be an abstract Wiener space and let
  $\{\pi_\theta,\,\theta\in[0,1]\}$ be a strictly increasing continuous
  resolution of the identity on $H$, and $\calF$ its induced filtration.
  If $Tw$ is a measure invariant transformation on $(W,H,\mu)$ then
  there exists a $\RRR\!\in\!\Lwa{W}$ which is a.s. an isometry on $H$,
  such that
 \begin{equation} \Label{4-5}
   Tw =\delta \RRR\, .
  \end{equation}
  Conversely if $\RRR\!\in\!\Lwa{W}$ is a.s. an isometry on $H$ then
  $\RRR\!\in\!\Domdel{2,W}$\, and $\delta\RRR$ is measure preserving.
 \end{theorem}
 \vspace{-.4cm}

 \noindent (Note that almost surely $R$'s range is contained in $H$, but
  its divergence is $W$--valued).\vspace{-.1cm}

\pf
 By our assumptions, every $ \eta_i $ can be uniquely represented as
 $ \eta_i =\delta u_i $ where the $ u_i $ are adapted, in the domain of
 $ \delta$, and $ u_i\in\DD_2 (H) $. Define $\boldR$ by
 \begin{equation} \label{4-6}
  \boldR(w) e_i = u_i
 \end{equation}
 then $ \boldR(w) $ is weakly adapted, and satisfies the assumptions of the
 previous result. Hence $\boldR$ is an isometry and
 $ Tw=\sum\,\delta(\boldR e_i) e_i $. In the converse direction, since
 $ \boldR(w) $ is weakly adapted, by Corollary~2.6.1 of \cite{uzbook},
 $ m_\theta =\delta (\pi_\theta \boldR h) \;\theta \in [0,1] $ is a
 $ \calF_\theta $ square integrable martingale and
 $ \lip m\rip_\theta = |\pi_\theta \boldR h|_H^2 $. Consequently by the Girsanov
 (or the stronger Novikov) condition
 \begin{align*}
  1 & = E\,\exp\left\{\delta(\boldR h) -\half |\boldR h|\right\} \\
%                                          \quad  \mathrm{and\ hence} \\
    & = E \, \exp\left\{\delta(\boldR h) - \half|h|^2\right\} \, .
 \end{align*}
 It follows that $ \delta (\boldR h) $ is $ N(0,|h|^2) $ and that
 $ \delta (\boldR e_i ) $ are i.i.d. $ N(0,1) $, so that
 \begin{equation*}
  Tw=\sum\,\delta (\RRR e_i)\,e_i = \delta \RRR \, .
 \end{equation*}
 \vspace{-1.7cm}

 \qed
 \section{Concluding Remarks}
  There is certainly no uniqueness in the representation of a random
  variable as a divergence if adaptedness of the integrand is not required.
  If a scalar random variable $\phi$, for example, can be written as\
  $\phi=\delta v$, and if
  \[ U_0=\{u\!\in\!\domdel{2},\ \delta u=0\} \]
  (that is, $U_0$ is the nonempty class of ``divergence free" integrands),
  then $\phi=\delta(v\!+\!u)$\ for any $u\!\in\!U_0$. The same is
  true for vector valued random variables.

  The question arises if there is a canonical integrand $\bar{v}$, for example
  \begin{equation} \Label{minimumenergy}
   E\|\bar{v}\|_{_H}^2=\min\left\{E|v|^2_{_H},\ \phi=\delta v\right\}
  \end{equation}
  or equivalently
  \[ E(v,u)_{_H}=0\ \ \ \forall u\!\in\!U_0\ \ \ \ \ \ \
                                \ ({\rm i.e.}\ v\!\in\!U_0^{\perp})\ .\]

  If we denote $L^2_{\rm e}(\mu;H):=\left\{\grad F,\
                             F\!\in\!\mathbb{D}_{2,1}\right\}$
  the space of {\em exact} $H$--valued random variables, then clearly
  $L^2_{\rm e}(\mu;H)\!\subset\!U_0^{\perp}$ since
  $E(\grad F,u)\!=\!EF\delta u$. Thus if $\phi=\delta(\grad F)$ for some
  $\grad F\!\in\!L^2_{\rm e}(\mu;H)$, then $\bar{v}=\grad F$ is the
  (necessarily unique) integrand which satisfies~\req{minimumenergy}.

  Let $\calL\!=\!\sum_{n=0}^\infty n\calP_n$ be the Ornstein--Uhlenbeck, or
  number, operator on $L^2(\mu)$, where $\calP_n$ is $L^2(\mu)$'s projection
  onto its $n$th homogeneous chaos, and ${\rm dom}\calL$ is the appropriate
  domain of convergence. From its definition, we see that $\calL$'s restriction
  to\ ${\rm dom}\calL\cap\{\phi\!\in\!L^2(\mu),\,E\phi=0\}$\ has a bounded
  inverse. In addition, it is well known that $\phi\!\in\!{\rm dom}\calL$\ \ if
  and only if\ $\phi\!\in\!\mathbb{D}_{2,1}$\ and\ $\grad\phi\in\domdel{}$,\ in
  which case $\calL\phi=\delta\grad\phi$.

  From the above discussion we conclude that
  \begin{equation} \Label{Ourep}
   \phi\!=\!E\phi\!+\!\delta\left(\grad\calL^{-1}(\phi-E\phi)\right),
  \end{equation}
  and that $\bar{v}\!=\!\grad\calL^{-1}(\phi-E\phi)$\ is the unique exact
  integrand in terms of which $\phi$ can be represented as a divergence, and as
  such satisfies the minimality condition~\req{minimumenergy}. Note that
  $\bar{v}$ is in general quite different from the adapted integrand discussed
  in this work; they coincide if and only if $\phi$ belongs to the first chaos\
  $\calP_1(L^2(\mu))$.

  The Ornstein--Uhlenbeck operator $\calL$ can be defined just as well in
  $L^2(\mu;B)$ for any Banach space $B$ (cf. for example \cite{S94}) via its
  interpretation as the generator of the Ornstein--Uhlenbeck semigroup. However,
  in order to extend~\req{Ourep} to $B$--valued $\phi$'s, assumptions on $B$
  seem to be needed in this case to conclude that $\calL$ has a bounded inverse
  on $L^2(\mu;B)$\,'s subspace of zero expectation, and this restricts the
  extension of the above argument when trying to obtain~\req{Ourep} for vector
  valued random variables.
%
%%%%%%%%%%%%%%%%%%%%%%%%%%%%%%%%%%%%%%%%%%%%%%%%%%%%%%%%%%%%%%%%%%%%%%%%%%%%
 
 \end{document}